\documentclass{article}
\usepackage[utf8]{inputenc}

\title{Alexander Quandles and Detecting Causality}
\author{Jack Leventhal}
\date{August 2022}
\usepackage{amsfonts}
\usepackage{times}
\usepackage{amssymb}
\usepackage{amsmath}
\usepackage{graphicx}
\usepackage{subcaption}
\usepackage{float}
\usepackage{array}
\newcolumntype{P}[1]{>{\centering\arraybackslash}p{#1}}
\newcolumntype{M}[1]{>{\centering\arraybackslash}m{#1}}
\begin{document}

\maketitle
\begin{abstract}
    In a recent paper, Allen and Swenberg investigated which link polynomials are capable of detecting causality in (2+1)-dimensional globally hyperbolic spacetimes. They ultimately suggested it is likely that the Jones Polynomial accomplishes this, while the Alexander-Conway polynomial, independently, is insufficient. As the Alexander-Conway polynomial on its own is not likely to detect causality an additional piece of information must be supplemented to the link polynomial. This paper aims to examine the ability of Alexander quandles, to distinguish the connected sum of two Hopf links and Allen-Swenberg Links. As the number of homomorphisms given by the Alexander quandles for the connected sum of Hopf links and the number of homomorphisms for the Allen-Swenberg Link are the same, we can conclude that the links are not distinguishable from one another using Alexander quandles and hence these quandles cannot capture causality when added to the Alexander-Conway polynomial.
\end{abstract}
\section{Introduction}
Mathematical \emph{knots} are embeddings of circles into Euclidean space (i.e. $S^1$ into $\mathbb{R}^3$), and can be continuously transformed via a series of Reidemeister moves. The most elementary of knots is a circle, commonly referred to as the \emph{trivial knot}, or \emph{unknot}. When a collection of knots is intertwined together, the result is a \emph{link}. By this definition, a knot can be considered as a link with a single component. Examples of links include the \emph{trivial link}, which is composed of two separate, unlinked circles and is also known as the \emph{unlink}. The simplest nontrivial link is the Hopf link, consisting of two circles linked together once.  
\newline
\newline
We can define a \emph{(2+1)-dimensional globally hyperbolic spacetime}, denoted $X$, as the result from conditioning the causal structure of a spacetime manifold, with a \emph{Cauchy surface} $\Sigma$ defined as a subset of spacetime which is crossed exactly once by every physically possible (causal) trajectory [HE73]. These Cauchy surfaces are defined as every causal curve, or a curve representing movement with speed less than or equal to the speed of light, which intersects the surface exactly once. Moreover, this equivalence between a globally hyperbolic spacetime and an existing Cauchy surface was proved by Geroch [Ger70] who furthermore showed that globally hyperbolic spacetimes are homeomorphic to $\Sigma \times \mathbb{R}$ (Bernal and Sanchez also proved globally hyperbolic spacetime to be diffeomorphic to $\Sigma \times \mathbb{R}$ [BS03]).
\newline
\newline
In order for global hyperbolicity to be maintained, Bernal and Sanchez laid out the two following conditions a spacetime must satisfy [BS07]:
\begin{enumerate}
    \item $\forall\ p, q \in X, J^+ (p) \cap J^- (q)$ is compact
    \item No time travel
\end{enumerate}
\noindent
The first condition states that for every two points $p$ and $q$ in globally hyperbolic spacetime $X$ the intersection of the causal future of point $p$ and the causal past of point $q$ must be compact; this requisite is also referred to as the \emph{absence of naked singularities}. The second condition requires our spacetime to lack the possibility of time travel. A \emph{causal trajectory} between two points $a$ and $b$ from a curve $\gamma$ connecting the two points must satisfy the inequality $\gamma '(t) \cdot \gamma '(t) < 0$, that is the dot product of the velocity vector of the curve with itself must be timelike or null, meaning that one can get from point $a$ to point $b$ without exceeding the speed of light. 
\newline
\newline
Utilizing the total space of spherical cotangent bundle $ST^{\ast}\Sigma$ of a Cauchy surface $\Sigma$, which in this case of a Cauchy surface being a plane will be homeomorphic to a solid torus (i.e. $S^1 \times \mathbb{R}^{2}$) we can identify it to the space of unparameterized future-directed null geodesics (light rays), $N$, in our globally hyperbolic spacetime, $X$. Doing so, we can subsequently associate the sphere of light rays passing through a point $x$ with the sphere $S_x \subset N_x$, called the \emph{sky} of $x$. For (2+1)-dimensional spacetimes the sky $S_x$ is homeomorphic to a circle and classified as a subset of the solid torus (i.e. $S^1 \times \mathbb{R}^{2}$). 
\newline
\newline
The notions of globally hyperbolic spacetimes and linking coalesce into the \emph{Low Conjecture}, which states that two events $x, y$ $\in X$ are causally related if and only if their skies $S_x \sqcup S_y$ are linked and that the resulting link is nontrivial. The Low Conjecture for (2+1)-dimensional spacetimes, where Cauchy surfaces $\Sigma$ are homeomorphic to $\mathbb{R}^{2}$, was proved by Chernov and Nemirovski [CN10], who demonstrated the relationship between causality and linking holds as long as the $\Sigma$ is not homeomorphic to $S^2$ or $\mathbb{R}P^2$. For higher dimensional spacetimes, Natário and Tod postulated the \emph{Legendrian Low Conjecture} that for (3+1)-dimensional spacetimes with $\Sigma$ diffeomorphic to $\mathbb{R}^{3}$, which states two events are causally related if and only if their skies intersect or are Legendrian linked [NT04]. This conjecture, comparable to the Low Conjecture though extended into higher dimensional spacetimes, was also proved by Chernov and Nemirovski [CN10]. The results of these conjectures beg the question of how we can confirm that the skies of two events ($S_x$ and $S_y$) are linked and therefore causally related.
\newline
\newline
Natário and Tod introduced a large number of families of skies which in pairs correspond to causally related events and are each associated with a link [NT04]. This link formed by a pair of skies correlates to a nontrivial Kauffman polynomial, which can be transformed into the Jones polynomial by a change in variables. 
\newline
\newline
From Chernov and Rudyak [CR08], the \emph{universal cover} $\tilde{X}$ of a (2+1)-dimensional globally hyperbolic spacetime $X$ with Cauchy surface $\Sigma$ is a globally hyperbolic spacetime with Cauchy surface $\tilde{\Sigma}$ as the universal cover of $\Sigma$. For two causally related points $x, y$ $\in X$, the lift of the path $\gamma$, $\tilde{\gamma}$, connects the lifts $\tilde{x}$, $\tilde{y}$ of $x, y$ to $X$, which can then be considered causally related in $\tilde{X}$. When the universal cover of $\Sigma$ is homeomorphic to $\mathbb{R}^{2}$ and is not $S^2$, the aforementioned universal cover, applied by Chernov, Martin, and Petkova [CMP20], shows that the Khovanov homology, a “categorification” of the Jones Polynomial [Bar02], can detect causality in (2+1)-dimensional spacetimes with $\Sigma \neq S^2$, $\mathbb{R}P^2$. Similarly, the Alexander-Conway polynomial, categorified by Heegaard Floer homology [KS16], was also proven by Chernov, Martin, and Petkova to detect causality in identical settings [CMP20]. The link polynomials described are strictly weaker link invariants than their respective homologies, though these polynomials may still detect causality. From the results of Allen and Swenberg [AS20], the Jones Polynomial is likely enough to detect causality, however, the Alexander-Conway polynomial likely cannot. In this paper, we explore whether the addition of the Alexander Quandle to the Alexander-Conway polynomial can capture causality. 
\section{Quandles}
A \emph{quandle} is defined as a set $X$ with a binary operation that satisfies the following three axioms:
\begin{enumerate}
    \item $x$ $\triangleright$ $x$ = $x$, for all $x \in X$
    \item For elements $x, y$ $\in X$, there exists some element $z$ such that $x$ = $y$  $\triangleright$ $z$
    \item $(x$ $\triangleright$ $y)$ $\triangleright$ $z$ = $(x$ $\triangleright$ $z)$ $\triangleright$ $(y$ $\triangleright$ $z)$, for all elements $x, y, z \in X$
\end{enumerate}
\noindent
The first axiom demonstrates the idempotency of quandles. The second alludes to the dual rack $\triangleright^{-1}$, the inverse operation such that $(x$ $\triangleright$ y) $\triangleright^{-1}$ $y$ $=$ $x$ for all elements $x, y$ $\in X$. One can also consider the second axiom by stating that for all $x$ in the set $X$ the map $f_x: X \to X$ is defined by $f_x (y)$ = $y$ $\triangleright$ $x$, and its inverse $f^{-1}_x (y)$, or $y$ $\triangleright^{-1}$ $x$. Finally, the third axiom exhibits the self-distributivity of the quandle. Besides these axioms it is worth noting that the quandle is non-commutative and non-associative and therefore generally $x$ $\triangleright$ $y$ $\neq$ $y$ $\triangleright$ $x$ and $(x$ $\triangleright$ $y)$ $\triangleright$ $z$ $\neq$ $x$ $\triangleright$ $(y$ $\triangleright$ $z)$ respectively. 
\newline 
\newline 
The relationship between knots and quandles is derived from the quandle relations at the crossings of a knot. Assigning variables to the arcs of the knot, which will ultimately be the elements of the quandle, allows relations to be made between arcs at each crossing, as depicted in the following diagrams:
\newpage
\begin{figure}[h]
\centering
\includegraphics[width=0.9\linewidth, height=5cm]{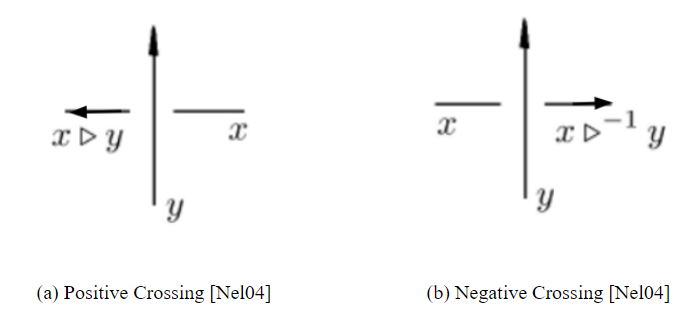}
\caption{Quandle Relations at Crossings}
\end{figure}
\noindent
The diagram on the left (a), represents the quandle coloring resulting from a positive crossing with the arc $x$ crossing under the arc labeled $y$ from right to left to form the strand colored by $x$ $\triangleright$ $y$. In a similar fashion, in (b) the arc $x$ crosses under the arc $y$ from left to right forming the strand colored by $x$ $\triangleright^{-1}$ $y$. These illustrations will later be applied to larger links with a substantial number of crossings and extensive quandle relationships.
\newline
\newline
Preliminarily, to ensure that quandle colorings are an invariant of knots, we must verify that quandle colorings remain unchanged after undergoing the three Reidemeister moves, that both diagrams prior and after the moves are equivalent. The axioms underlying the quandle translate to the Reidemeister moves as follows:
\begin{figure}[h]
\centering
\includegraphics[width=0.9\linewidth, height=5cm]{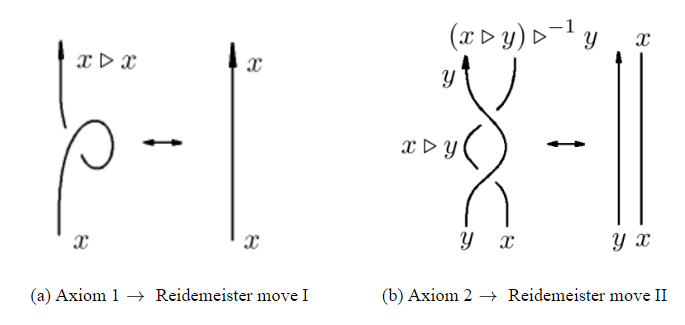}
\end{figure}

\newpage
\begin{figure}[h]
\centering
\includegraphics[width=0.8\linewidth, height=6cm]{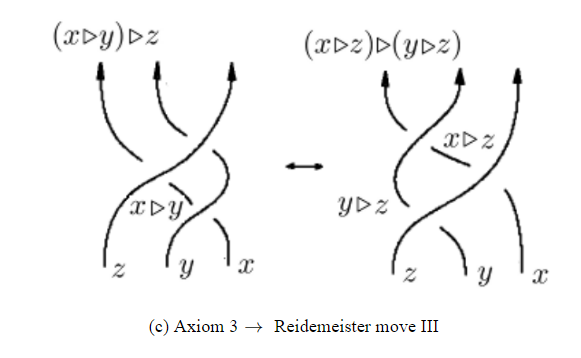}
\caption{Reidemeister Moves and Corresponding Quandle Axioms}
\end{figure}

\noindent From the three figures above, (a) exemplifies how the arc on the left, $x$ $\triangleright$ $x$, is equivalent to $x$ through the first Reidemeister move, additionally stated as the first axiom of the quandle. The second figure, (b), depict the right inverse with the arcs $(x$ $\triangleright$ $y)$ $\triangleright^{-1}$ $y$ on the left and the $x$ strand on the right being equivalent as previously mentioned in the second axiom. Lastly, the third figure (c) demonstrates the self-distributivity for quandles with equivalent arcs $(x$ $\triangleright$ $y)$ $\triangleright$ $z$ and $(x$ $\triangleright$ $z)$ $\triangleright$ $(y$ $\triangleright$ $z)$.
\newline
\newline
Therefore, if labeled arcs of a quandle are constructed according to the crossing structure, another labeling after a Reidemeister move corresponds to the same fundamental quandle. Furthermore, we can gather that the number of ways to label a knot diagram, by elements satisfying the labeling conditions, will be the same for two diagrams of the knot. Consequently, in order to distinguish one knot from another, the number of ways to label the two knots' diagrams by elements of a fixed, finite quandle, must be different. If the number is equivalent for both diagrams, then no distinction can be made. This invariant, known as the \emph{quandle counting invariant}, is denoted $|\emph{Hom} (Q(K), Q)|$, for some knot diagram $K$. 

\section{Useful Quandles \& Invariants}

Applying a quandle to a knot or link we obtain a system of equations based off of the elements of the quandle, the arcs enumerated as variables in a knot or link. Solving the system of equations in some modular class, we acquire the number of \emph{homomorphisms} given by the fundamental quandle, simply the number of solutions to our system. Quandle homomorphisms can be classified by the property $f(x$ $\triangleright$ $y)$ = $f(x)$ $\triangleright$ $f(y)$. Additionally, given a quandle $X$, and a link $L$, we can define a \emph{coloring} of $L$ by $X$ to be a replacement of the elements in $X$ instead of the variables in $Q(L)$ such that relationships that comprise the fundamental quandle are maintained. This can also be considered to be a homomorphism from $Q(L)$ into $X$. Once more it is precisely the number of homomorphisms that will allow us to distinguish between links.
\newline
\newline
One example of a quandle is the \emph{Takasaki quandle} defined by the operation $x$ $\triangleright$ $y$ = $2y - x$. Checking the axioms of quandles we can see the following:
\begin{enumerate}
    \item $x$ $\triangleright$ $x$ = $2x - x$ = $x$, therefore the Takasaki quandle satisfies the first condition
    \item For every $x, y$ in the set of quandles $X$ there will be some element $z$ $\in X$, such that $x$ = $y$ $\triangleright$ $z$ = $2z - y$ and $(x$ $\triangleright$ $y)$ $\triangleright^{-1}$ $y$ = $x$
    \item $(x$ $\triangleright$ $y)$ $\triangleright$ $z$ = $2z - (x$ $\triangleright$ $y)$ = $2z - (2y - x)$ = $x - 2y + 2z$ and \newline $(x$ $\triangleright$ $z)$ $\triangleright$ $(y$ $\triangleright$ $z) = 2(y$ $\triangleright$ $z) - (x$ $\triangleright$ $z)$ $=$ $2(2z – y) - (2z – x)$ $=$ $x – 2y + 2z$
    \newline
    \newline
    As the above two expressions are equivalent, we have shown that the self-distributivity law holds for the Takasaki quandle.
    
\end{enumerate}
\noindent
In addition to the Takasaki quandle, the \emph{Alexander quandle} is derived from \emph{(t, s)-rack structures}. \emph{Racks} are more generalized versions of quandles which only satisfy the second and third of the previously described axioms. As laid out by Nelson [Nel14] the (t, s)-rack structure can be defined as follows: for $\Lambda$, the ring $\mathbb{Z}$[t, $t^{-1}$, s] modulo the ideal generated by $s^2 - (1 - t)s$. We then define the operation $x$ $\triangleright$ $y = tx + sy$, and take M $= \mathbb{Z}_{n}$ with $t, s \in M$ such that $gcd(n, t) = 1$ and $s^2 = (1 - t)s$. If $s = 1 - t$, then M is an Alexander quandle. From this substitution we could also write the operation for the Alexander quandle as $x$ $\triangleright$ $y = tx + (1 - t)y$. By similar logic used for the Takasaki quandle, we can also see how the Alexander quandle satisfies the three axioms:
\begin{enumerate}
    \item For all $x \in X$, $x$ $\triangleright$ $x = tx - (1 - t)x = x$
    \item For every $x, y \in X$ there will be some element $z \in X$, such that $x = y $ $\triangleright$ $z = ty + (1 - t)z$ and $(x$ $\triangleright$ $y)$ and $(x$ $\triangleright$ $y)$ $\triangleright^{-1}$ $y = x$
    \item A straightforward computation additionally confirms that for $x, y, z \in X$, $(x$ $\triangleright$ $y)$ $\triangleright$ $z = (x$ $\triangleright$ $z)$ $\triangleright$ $(y$ $\triangleright$ $z)$
\end{enumerate}
\noindent
In order to produce stronger invariants that will allow us to further distinguish knots and links, we can construct an \emph{enhanced linking polynomial} from the colorings. For each coloring $f \in Hom(Q(L), X)$ in quandle $X$ and link $L$, we count the number of elements that appear in the coloring, denoted $|Im(f)|$. Now we can define the \emph{enhanced linking polynomial}, in variable q, with the formula $\phi_x (L) = \sum_{f \in Hom(Q(L), X)} q^{|Im(f)|}$. The \emph{enhanced linking polynomial} is an invariant of L and can be particularly useful in distinguishing links.
\newline
\newline
For example, we can complete the process of finding the fundamental quandle, Alexander quandle, number of colorings, and enhanced linking polynomial for the trefoil, a relatively simple nontrivial knot.  
\newline
\newline
Below we can add orientation to a right-handed trefoil and label each arc and crossing.  
\newline
\newline
\newpage
\begin{figure}[h]
\centering
\includegraphics[width=0.6\linewidth, height=6.5cm]{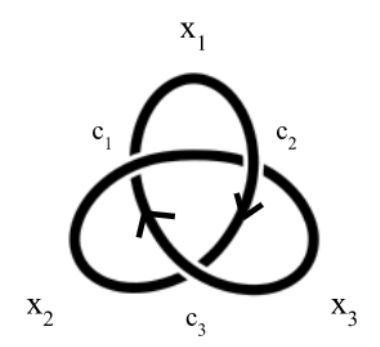}
\caption{Trefoil Knot}
\end{figure}

\noindent
This instance we will use the Alexander quandle derived from a (2, 2)-rack structure (i.e. with $t$ and $s$ equal to 2), and consider the set $\mathbb{Z}_{3}$. Notice this is equivalent to an Alexander quandle with $n = 3$ and $t = 2$ as $gcd(3, 2) = 1, s^2 = 4 \equiv 1$ $(mod$ $3)$ is equivalent to $(1 - t)s = (1 - 2)2 = -2 \equiv 1$ $(mod$ $3)$, and $s = 2 \equiv 1 - t = 1 - 2 = -1 \equiv 2$ $(mod$ $3)$.
\newline
\newline 
In the table below, we will show the quandle relation at each crossing in addition to the associated Alexander quandle, which here is the operation $x \triangleright y = 2x + (1 - 2)y = 2x - y \equiv 2x + 2y$ $(mod$ $3)$.
\newline
\begin{center}
\renewcommand{\arraystretch}{1.75}
\begin{tabular}{|>{\centering\arraybackslash}p{1.75cm}|wc{4cm}|wc{5cm}|}
  \hline
  Crossing Number & Fundamental Quandle & Alexander Quandle \\ 
  \hline
  $c_1$ & $x_3 = x_1 \triangleright x_2$ & $x_3 = 2x_1 + 2x_2$ $(mod$ $3)$  \\ 
  \hline
  $c_2$ & $x_2 = x_3 \triangleright x_1$ & $x_2 = 2x_3 + 2x_1$ $(mod$ $3)$ \\ 
  \hline
  $c_3$ & $x_1 = x_2 \triangleright x_3$ & $x_1 = 2x_2 + 2x_3$ $(mod$ $3)$ \\
  \hline
\end{tabular}
\end{center}
\noindent
\newline
\noindent
Rearranging any of the equations for the Alexander Quandle the result is identical, leaving $x_1 + x_2 + x_3 = 0$ $(mod$ $3)$. This equation will determine the number of homomorphisms and colorings for the trefoil knot. Quickly computing the solutions, we get $\{x_1, x_2, x_3\} = \{\{0, 0, 0\}, \{1, 1, 1\}, \{2, 2, 2\}, \{0, 1, 2\}, \{0, 2, 1\}, \{1, 0, 2\}, \{1, 2, 0\}, \{2, \newline 0, 1\}, \{2, 1, 0\}\}$, for a total of 9 colorings. Looking at the solutions more closely three solutions are trivial with all variables equivalent to one another; they can also be considered monochromatic (i.e. $|Im(f)| = 1$). We also have six other colorings with all variables being different values (or rather different colors), called tri-colorings as three different values, or colors, were used in the solution set. From the formula for the enhanced linking polynomial invariant, we will have one term $3q$ representing the three monochromatic colorings and another term $6q^3$ for the six tricolors. Therefore, the $\Phi - polynomial$ for the trefoil ($3_1$) is $\Phi_{\mathbb{Z}_{3}} (3_1) = 3q + 6q^3$.

\section{Alexander Quandles \& Detecting Causality}

The connected sum of two Hopf links originates from the knots formed within solid tori ($S^1 \times \mathbb{R}^{2}$) in (2+1)-dimensional globally hyperbolic spacetimes. When the Cauchy surface is homeomorphic to $\mathbb{R}^{2}$, the skies of two causally unrelated events will be circles in the solid torus, isotopic to the unlinked pair of longitudes. Using these three-component links (two from the skies, the third from deleting the trivial knot that makes $\mathbb{R}^{3}$ a doughnut), we study the connected sum of Hopf links formed in this situation to see what can distinguish this link, for our purposes from Allen-Swenberg links.  
\newline
\newline
Two examples from Allen and Swenberg vital in studying causality are the connected sum of Hopf links characterized previously, and the Allen-Swenberg Link, both of which are below. After finding the Alexander quandles for both links we will compare the number of homomorphisms and colorings for the two and examine their enhanced linking polynomials.  
\newline

\begin{figure}[h]
\centering
\includegraphics[width=0.9\linewidth, height=5cm]{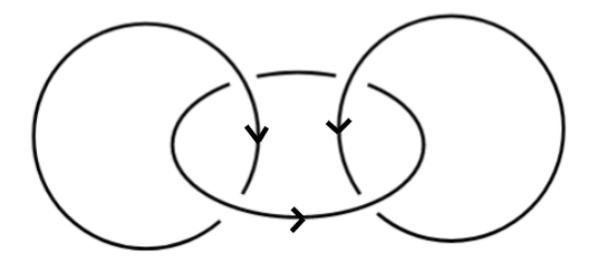}
\caption{Connected sum of Hopf links [AS20]}
\end{figure}

\newpage
\begin{figure}[h]
\centering
\includegraphics[width=0.9\linewidth, height=12cm]{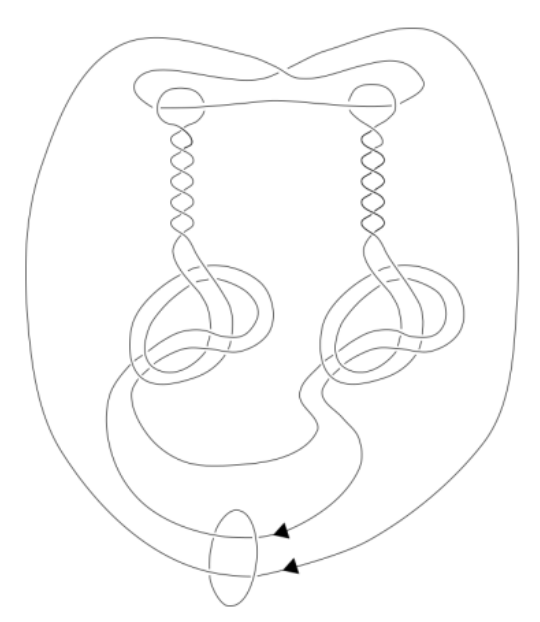}
\caption{Allen-Swenberg Link [AS20]}
\end{figure}

\noindent The Alexander quandle used can be constructed from a (t, s)-rack structure with the substitution $s = 1 - t$, all over $\mathbb{Z}_{n}$. Therefore we have the quandle relation $x$ $\triangleright$ $y = tx + (1 - t)y$. 
\newline
\newline
Similar to the process for the trefoil, we will label the arcs and crossing numbers, first for the connected sum of Hopf links, as shown in the diagram below.  

\newpage
\begin{figure}[h]
\centering
\includegraphics[width=1\linewidth, height=5cm]{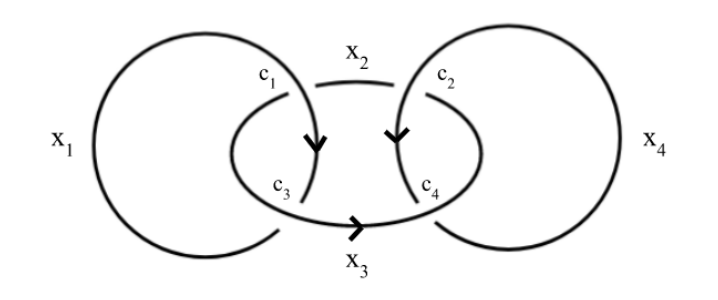}
\caption{Labeled diagram of the connected sum of Hopf links}
\end{figure}

\noindent The table hereunder shows for each crossing the element of the fundamental quandle as well as the corresponding equation from the Alexander quandle operation. 
\newline

 \begin{center}
\renewcommand{\arraystretch}{1.75}
\begin{tabular}{|>{\centering\arraybackslash}p{3cm}|wc{4cm}|}
  \hline
  Crossing Number & Fundamental Quandle \\ 
  \hline
  $c_1$ & $x_2 = x_3 \triangleright x_1$ \\ 
  \hline
  $c_2$ & $x_3 = x_2 \triangleright x_4$ \\ 
  \hline
  $c_3$ & $x_1 = x_1 \triangleright x_3$ \\
  \hline
  $c_4$ & $x_4 = x_4 \triangleright x_3$ \\
  \hline
\end{tabular}
\end{center}
\noindent
\newline
The corresponding (non-involutory) Alexander quandle for each element of the fundamental quandle relation of the connected sum of Hopf links is as follows:
\begin{center}
$
\begin{cases} x_2 = t \cdot x_3 + (1 - t) \cdot x_1 \\ x_3 = t \cdot x_2 + (1 - t) \cdot x_4 \\ x_1 = t \cdot x_1 + (1 - t) \cdot x_3 \\ x_4 = t \cdot x_4 + (1 - t) \cdot x_3 \end{cases}  
$
\end{center}
This can subsequently manipulated to ensure the equations are in standard form as follows:
\begin{center}
$
\begin{cases} t \cdot x_3 + (1 - t) \cdot x_1 - x_2 = 0 \\ t \cdot x_2 + (1 - t) \cdot x_4 - x_3 = 0 \\ (t - 1) \cdot x_1 + (1 - t) \cdot x_3 = 0 \\ (t - 1) \cdot x_4 + (1 - t) \cdot x_3 = 0 \end{cases}  
$
\end{center}
This corresponds to the following matrix:
\newline
\newline
\begin{center}
$
\begin{bmatrix}
1 - t & -1 & t & 0 & 0\\
0 & t & -1 & 1 - t & 0\\
t - 1 & 0 & 1 - t & 0 & 0\\
0 & 0 & 1 - t & t - 1 & 0\\
\end{bmatrix}
$
\end{center}
Which ultimately reduces to the following:
\newline
\begin{center}
$
\begin{bmatrix}
1 & 0 & 0 & -1 & 0\\
0 & 1 & 0 & -1 & 0\\
0 & 0 & 1 & -1 & 0\\
0 & 0 & 0 & 0 & 0\\
\end{bmatrix}
$
\end{center}
Hence, all variables are shown to be equivalent, and therefore the sole colorings for the connected sum of Hopf links are trivial. Furthermore, as the possible values for each variable ranges from $0$ to $n - 1$ over $\mathbb{Z}_{n}$, the enhanced linking polynomial for the sum of two Hopf links under the Alexander quandle is $\Phi_{\mathbb{Z}_{n}} (2^2_1$ $\#$ $2^2_1) = nq$, meaning that there are $n$ possible colorings of the link over $\mathbb{Z}_{n}$.
\newline
\newline
We now repeat this process for the Allen-Swenberg link, labeling the arcs and crossing numbers, adding orientation to the link. The same Alexander quandle relation given by $x$ $\triangleright$ $y = tx + (1 - t)y$ will be utilized for this computation. 

\newpage
\begin{figure}[h]
\centering
\includegraphics[width=1\linewidth, height=13cm]{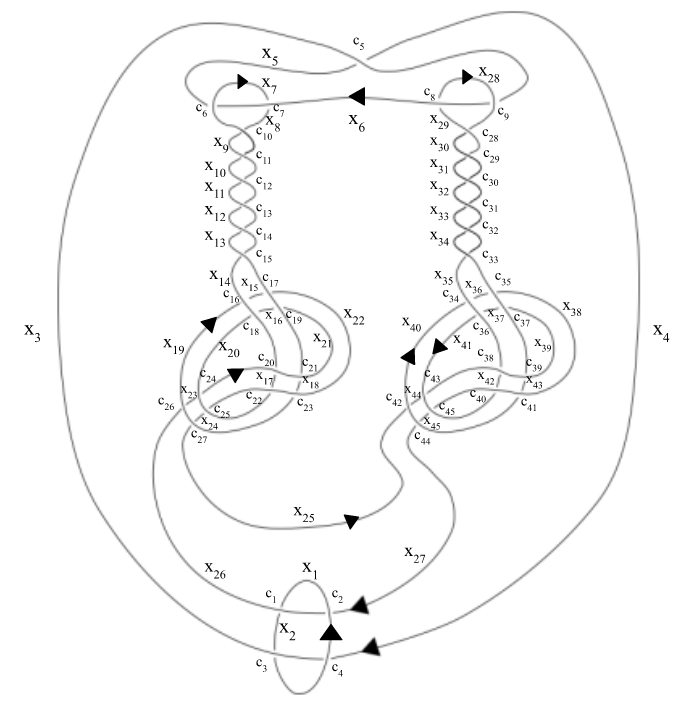}
\caption{Labeled Allen-Swenberg Link}
\end{figure}

\noindent As with previous examples the table records the element of the fundamental quandle and the corresponding Alexander quandle at each crossing of the link.

\begin{center}
\renewcommand{\arraystretch}{1.75}
\begin{tabular}{|>{\centering\arraybackslash}p{1.75cm}|wc{4cm}|}
  \hline
  Crossing Number & Fundamental Quandle\\ 
  \hline
  $c_1$ & $x_2 = x_1 \triangleright x_{26}$\\ 
  \hline
  $c_2$ & $x_{26} = x_{27} \triangleright x_1$\\ 
  \hline
  $c_3$ & $x_1 = x_2 \triangleright x_3$ \\
  \hline
  $c_4$ & $x_3 = x_4 \triangleright x_1$  \\
  \hline
  $c_5$ & $x_4 = x_5 \triangleright x_3$ \\ 
  \hline
  $c_6$ & $x_5 = x_6 \triangleright x_7$ \\ 
  \hline
  $c_7$ & $x_8 = x_7 \triangleright x_6$ \\ 
  \hline
  $c_8$ & $x_{29} = x_{28} \triangleright x_6$ \\
  \hline
  $c_9$ & $x_3 = x_6 \triangleright x_{28}$ \\ 
  \hline
  $c_{10}$ & $x_9 = x_8 \triangleright x_7$ \\ 
  \hline
  $c_{11}$ & $x_7 = x_{10} \triangleright x_9$ \\ 
  \hline
  $c_{12}$ & $x_{11} = x_9 \triangleright x_{10}$ \\
  \hline
  $c_{13}$ & $x_{10} = x_{12} \triangleright x_{11}$ \\
  \hline
  $c_{14}$ & $x_{13} = x_{11} \triangleright x_{12}$ \\
  \hline
  $c_{15}$ & $x_{12} = x_{14} \triangleright x_{13}$ \\
  \hline
  $c_{16}$ & $x_{19} = x_{15} \triangleright x_{14}$ \\
  \hline
  $c_{17}$ & $x_{22} = x_{15} \triangleright x_{13}$ \\
  \hline
  $c_{18}$ & $x_{20} = x_{16} \triangleright x_{14}$ \\
  \hline
  $c_{19}$ & $x_{21} = x_{16} \triangleright x_{13}$ \\
  \hline
  $c_{20}$ & $x_{14} = x_{17} \triangleright x_{21}$ \\
  \hline
  $c_{21}$ & $x_{13} = x_{18} \triangleright x_{21}$ \\ 
  \hline
  $c_{22}$ & $x_{20} = x_{17} \triangleright x_{22}$ \\ 
  \hline
    \end{tabular}
\end{center}
\begin{center}
\renewcommand{\arraystretch}{1.75}
\begin{tabular}{|>{\centering\arraybackslash}p{1.75cm}|wc{4cm}|}
  \hline
  $c_{23}$ & $x_{19} = x_{18} \triangleright x_{22}$  \\ 
  \hline
  $c_{24}$ & $x_{21} = x_{23} \triangleright x_{20}$  \\ 
  \hline
  $c_{25}$ & $x_{22} = x_{24} \triangleright x_{20}$ \\ 
  \hline
  $c_{26}$ & $x_{26} = x_{23} \triangleright x_{19}$ \\ 
  \hline
  $c_{27}$ & $x_{25} = x_{24} \triangleright x_{19}$ \\ 
  \hline
  $c_{28}$ & $x_{30} = x_{28} \triangleright x_{29}$ \\ 
  \hline
  $c_{29}$ & $x_{29} = x_{31} \triangleright x_{30}$ \\ 
  \hline
  $c_{30}$ & $x_{32} = x_{30} \triangleright x_{31}$ \\ 
  \hline
  $c_{31}$ & $x_{31} = x_{33} \triangleright x_{32}$ \\ 
  \hline
  $c_{32}$ & $x_{34} = x_{32} \triangleright x_{33}$ \\ 
  \hline
  $c_{33}$ & $x_{33} = x_{35} \triangleright x_{34}$ \\ 
  \hline
  $c_{34}$ & $x_{40} = x_{36} \triangleright x_{35}$ \\ 
  \hline
  $c_{35}$ & $x_{38} = x_{36} \triangleright x_{34}$ \\ 
  \hline
   $c_{36}$ & $x_{41} = x_{37} \triangleright x_{35}$ \\ 
  \hline
  $c_{37}$ & $x_{39} = x_{37} \triangleright x_{34}$ \\ 
  \hline
  $c_{38}$ & $x_{35} = x_{42} \triangleright x_{39}$ \\ 
  \hline
  $c_{39}$ & $x_{34} = x_{43} \triangleright x_{39}$  \\ 
  \hline
  $c_{40}$ & $x_{41} = x_{42} \triangleright x_{38}$\\ 
  \hline
  $c_{41}$ & $x_{40} = x_{43} \triangleright x_{38}$ \\ 
  \hline
  $c_{42}$ & $x_{25} = x_{44} \triangleright x_{40}$ \\ 
  \hline
  $c_{43}$ & $x_{39} = x_{44} \triangleright x_{41}$ \\ 
  \hline
  $c_{44}$ & $x_{27} = x_{45} \triangleright x_{40}$  \\ 
  \hline
  $c_{45}$ & $x_{38} = x_{45} \triangleright x_{41}$\\ 
  \hline
 \end{tabular}
 \end{center}
Similarly we can perform this process for the system of equations corresponding to the revised non-involutory Alexander quandle of the Allen-Swenberg Link. As there are $45$ equations involved the reduced row echelon form of the matrix can be calculated using Wolfram Mathematica [Wol22]. 
\begin{figure}[h]
\centering
\includegraphics[width=0.75\linewidth, height=7.5cm]{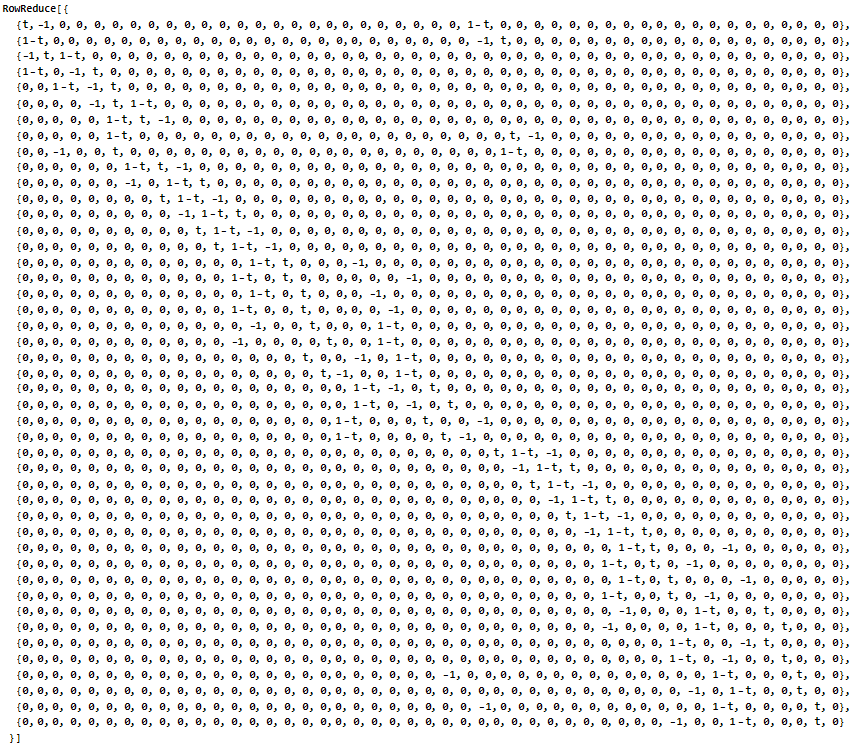}
\caption{Matrix for Alexander Quandle of Allen-Swenberg Link}
\end{figure}
\newline 
This matrix ultimately row reduces to a matrix similar to the form of the reduced row echelon form of the matrix for the connected sum of Hopf links with the elements of the main diagonal equal to $1$ (with the exception of the last row consisting entirely of zeroes) and the $45$th column (again with the exception of the last row) having all elements equal to $-1$. The $45$ $\times$ $46$ matrix is the following:
\newpage
\begin{figure}[h]
\centering
\includegraphics[width=0.75\linewidth, height=7.5cm]{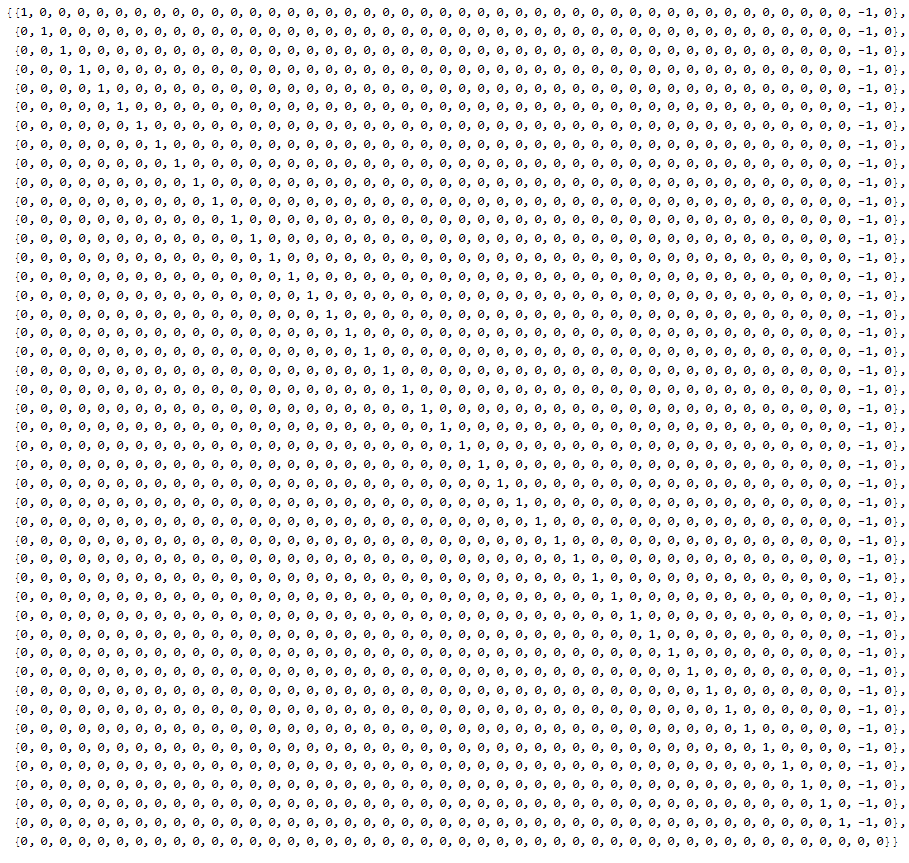}
\caption{RREF of Matrix for Alexander Quandle of Allen-Swenberg Link}
\end{figure}
\noindent
This solution demonstrates that the sole solutions from this system are trivial, with all variables equivalent to one another (i.e. $x_1 = x_2 = ... = x_{45}$). Hence, the enhanced linking polynomial for the Allen-Swenberg link over $\mathbb{Z}_{n}$ will be $nq$. As both links have the same number of trivial solutions and the same enhanced linking polynomial the two cannot be distinguished solely by the Alexander quandle.
\newline
\newline
We now examine the ability of the involutory quandle, when added to the Alexander quandle, to distinguish the two links. The involutory quandle equates the quandle operation with its inverse, or $\triangleright = \triangleright^{-1}$. Therefore the identity $(x \triangleright y) \triangleright^{-1} y = x$ is equivalent to $(x \triangleright y) \triangleright y = x$. Applying this to Alexander quandles, substituting the value of $x \triangleright y$ for $x$, or equivalently substituting $tx + (1 - t)y$ for x, by the definition of the Alexander quandle, we get the following:
\begin{center}
$t(tx + (1 - t)y) + (1 - t)y = x$
\end{center}
\begin{center}
$\Rightarrow$ $x \cdot t^2 - y \cdot t^2 + y = x$
\end{center}
\begin{center}
$\Rightarrow$ $(x - y) \cdot (t^2 - 1) = 0$
\end{center}
Therefore $x = y$ or $t = \pm 1$. We can now consider each of these cases. \newline
\newline
Case 1: $x = y$ 
\newline
If $x = y$ we can make the substitution into the Alexander quandle operation with $x$ $\triangleright$ $y = x$ $\triangleright$ $x = x = z$. Therefore we have shown that all elements would be equivalent to one another, only producing trivial solutions for both links, unable to distinguish them. 
\newline
\newline
Case 2: $t = 1$
\newline
If $t = 1$ we can make this substitution in the Alexander quandle operation to get $x \triangleright y = 1x + (1 - 1)y = x + 0y = x = z$. Therefore, regardless of $y$, for each quandle relation of the form $x \triangleright y = z$, $x = z$. In the connected sum of Hopf links this implies that $x_2 = x_3$, while there are no such equivalences for $x_1$ and $x_4$ (besides being equivalent to themselves). Therefore there are three free variables that can be used to describe the number of solutions to this system (i.e. one for $x_1$, another for $x_4$, and the last equivalent to $x_2$ and $x_3$) and therefore over $\mathbb{Z}_{n}$ there exists $n^3$ solutions (as there are n choices for the value of each free variable ranging from $0$ to $n - 1$). Similarly, the fact that $x = z$ can be used to show that for the Allen-Swenberg link the following equivalences hold:
\begin{enumerate}
    \item $x_1 = x_2$,
    \item $x_3 = x_4 = x_5 = x_6$, and
    \item $x_7 = x_8 = ... = x_{45}$
\end{enumerate}
With these three equivalences, over $\mathbb{Z}_{n}$ we have $n^3$ solutions, the same number of solutions as that of the connected sum of Hopf links. Hence, the number of homomorphisms of the connected sum of Hopf links and the Allen-Swenberg Link are the same and the two cannot be distinguished in this case.
\newline
\newline
Case 3: $t = -1$
\newline
As $t \in \mathbb{Z}_{n}$, we can use the value of $t = n - 1$ in our calculations. Under the involutory condition for our Alexander quandle relation this means the following equation must be true:
\begin{center}
    $(x \triangleright y) \triangleright y = (n - 1) \cdot ((n - 1) \cdot x + (1 - (n - 1)) \cdot y) + (1 - (n - 1)) \cdot y = x$
\end{center}
\begin{center}
    $\Rightarrow$ $(n - 1)^2 \cdot x + (n - 1) \cdot (2 - n) \cdot y + (2 - n) \cdot y = x$
\end{center}
\begin{center}
    $\Rightarrow$ $-y \cdot n^2 + 2y \cdot n + x \cdot n^2 - 2x \cdot n + x = x$
\end{center}
\begin{center}
    $\Rightarrow$ $n^2 \cdot (x - y) - 2n \cdot (x - y) = 0$
\end{center}
\begin{center}
    $\Rightarrow$ $n \cdot (n - 2) \cdot (x - y) = 0$
\end{center}
From the above equation we have three subcases now to consider:
\begin{enumerate}
    \item Firstly, n cannot equal zero as the remainder class $\mathbb{Z}_{0}$ is undefined. 
    \item If n = 2, this would correspond to the case in which t = 1, which we have proved cannot distinguish the links.
    \item Finally, if $x = y$, as shown in the first case, all elements of the quandle are equal and the solution solely consists of trivial solutions. 
\end{enumerate}
 Examining each of the possible cases under which we can add the involutive condition to the Alexander quandle, we have shown that the Alexander quandle when paired with the involutory quandle produces the same number of homomorphisms for both links and therefore cannot distinguish them. 
 \newline
 \newline
 Hence, neither the Alexander quandle nor the addition of the involutory quandle to the Alexander quandle can distinguish the connected sum of Hopf links and the Allen-Swenberg link.

\section{Conclusion}
From the results of this study, both the number of colorings, homomorphisms, and enhanced linking invariant for the connected sum of Hopf links and the Allen-Swenberg links will be the exact same with the addition of Alexander quandles and when the involutory quandle is supplemented to the Alexander quandle. Therefore, neither the Alexander quandles independently, nor the Alexander and involutory quandles together, when paired with the Alexander-Conway polynomial, are enough to detect causality in this scenario.   

\section{Acknowledgements}
This project was completed in the Summer of 2022 as a part of the Horizon Academic Research Program. The program was supervised under Professors Vladimir Chernov of Dartmouth College and Emanuele Zappala of Yale University. I would like to thank Professor Chernov and Dr. Zappala for their invaluable support as research mentors.  

\section{References}
[AS20] S. Allen, J. Swenberg, “Do Link Polynomials Detect Causality In Globally Hyperbolic Spacetimes?”, \emph{J. Math. Phys.} Vol. 62, No.3, 032503 (2021) 
\newline
\newline
[BS03] A. Bernal, M. Sanchez, “On smooth Cauchy hypersurfaces and Geroch's splitting theorem”, \emph{Commun. Math. Phys.} 243 (2003) 461 - 470  
\newline
\newline
[BS07] A. Bernal, M. Sanchez, “Globally hyperbolic spacetimes can be defined as "causal" instead of "strongly causal"”, \emph{Class. Quant. Grav.} 24 (2007) 745 – 750 
\newline
\newline
[CMP20] V. Chernov, G. Martin, I. Petkova, “Khovanov homology and causality in spacetimes”, \emph{J. Math. Phys.} 61, 022503 (2020) 
\newline
\newline
[CN10] V. Chernov, S. Nemirovski, “Legendrian links, causality, and the Low conjecture”, \emph{Geom. Funct. Anal.} 19 (2010), 1323 – 1333 
\newline
\newline
[CR08] V. Chernov, Y. Rudyak, “Linking and causality in globally hyperbolic spacetimes”, \emph{Commun. Math. Phys.} 279: 309 – 354, 2008 
\newline
\newline
[Ger70] R.Geroch, “Domain of dependence”, \emph{J. Math. Phys.} 11 (1970) 437 – 449 
\newline
\newline
[HE73] S. Hawking, G. Ellis, “The large scale structure of space-time", \emph{Camb. Mono. Math. Phys.} 1, Cambridge University Press (1973) 
\newline
\newline
[KS16] L. Kauffman, M. Silvero, “Alexander-Conway Polynomial State Model and Link Homology”, \emph{J. Knot Theo. Ram.} 25.3 (2016) 
\newline
\newline
[Kho99] M. Khovanov, “A categorification of the Jones Polynomial”, \emph{Duke J. Math}, 101 (2000), No. 3, 359 - 426 
\newline 
\newline
[NT04] J. Natário, P. Tod, “Linking, Legendrian linking and causality”, \emph{Proc. Lond. Math. Soc.} 88 (2004) 251 – 272 
\newline
\newline
[Nel04] S. Nelson, "Quandles and Racks", \newline https://www1.cmc.edu/pages/faculty/VNelson/quandles.html. (2004) 
\newline
\newline
[Nel14] S. Nelson, “Link invariants from finite racks”, \emph{Fund. Math.} 225 (2014) 234 – 258 
\newline
\newline
[Wol22] Wolfram Research, Inc. (www.wolfram.com), \emph{Mathematica Online, Version 13.1.} Champaign, IL (2022) 

\end{document}